\documentclass[12pt]{amsart}
\usepackage{latexsym,amsmath,amssymb,amscd}
\usepackage{eucal}
\theoremstyle{plain}
\newtheorem*{theorem}{Theorem}
\newtheorem*{corollary}{Corollary}
\newtheorem*{lemma}{Lemma}

\theoremstyle{definition}

\theoremstyle{remark}

\newcommand{\bbc}{\mathbb{C}}

\newcommand{\fb}{\mathfrak{B}}
\newcommand{\fm}{\mathfrak{M}}

\newcommand{\ra}{\rightarrow}
\newcommand{\norm}[1]{\| #1 \| }

\bibliography{/Users/jb/Praca/Styles/biblioteka.bib}
\usepackage{pb-diagram}
\usepackage{eucal}
\begin{document}
\title[]{Rapid decay and the metric approximation property }
\author{Jacek Brodzki}
\address{School of Mathematical Sciences, University of Southampton, Southampton SO17 1BJ}
\email{j.brodzki@soton.ac.uk, g.a.niblo@soton.ac.uk.}
\author{Graham A. Niblo}
\maketitle

A $C^*$-algebra $A$ is said to have the Metric Approximation Property (MAP) if the identity map $\text{id}: A\ra A$ can be approximated in the point-norm topology by a net of  finite rank contractions.
 The main purpose of this note is to prove the following theorem. 
 \begin{theorem}[]
 Let $\Gamma$ be a discrete group satisfying the rapid decay property with respect to a length function $\ell$
 which is conditionally negative. Then the reduced $C^*$-algebra $C^*_r(\Gamma)$ has the metric 
 approximation property. 
 \end{theorem}
 The central point of our proof is an observation that the proof of the same property for free groups due
 to Haagerup
 \cite{Haagerup:Inventiones} transfers directly to this more general situation. We also note that 
 under the same hypotheses, the Fourier algebra $A(\Gamma)$ has a bounded approximate identity, which implies that it too has the MAP. 
 
 A discrete group $\Gamma$ satisfies property (RD) (Rapid Decay) with respect to 
 a length function $\ell$ on $\Gamma$  if the operator norm of any element of the group ring 
 can be uniformly majorised by a Sobolev norm determined by $\ell$. In detail, this means the following. 
 
 The left action of a group $\Gamma$ on itself extends to the convolution action of the group 
 ring $\bbc \Gamma$  on the Hilbert space $\ell^2(\Gamma)$. This is the left regular representation $\lambda$ of $\Gamma$ which embeds the group ring in the $C^*$-algebra $\fb(\ell^2(\Gamma))$ of all bounded linear operators on $\ell^2(\Gamma)$.  The reduced 
 $C^*$-algebra $C^*_r(\Gamma)$ is the $C^*$-subalgebra of  $\fb(\ell^2(\Gamma))$ 
 generated by $\lambda (\bbc \Gamma)$.

 
 For any positive real number $s$ we define a Sobolev norm associated with the length function 
 $\ell$ by 
 $$
 \norm{f}_{\ell,s} = \sqrt{\sum_{\gamma\in \Gamma} |f(\gamma)|^2 (1+\ell(\gamma))^{2s}}
 $$
 for any $f\in\bbc \Gamma$. 
  
 Following Jolissaint \cite{Jolissaint} (see also \cite{CR}) we say that $\Gamma$ has property 
 (RD) with respect to the length function $\ell$ if and only if it satisfies the following property: There exist a $C> 0$ and $s> 0$ such that for all $f\in \bbc \Gamma$ 
 $$
 \norm{\lambda (f)} \leq C\norm{f}_{\ell, s}. 
 $$
where the norm on the left hand side is the operator norm in $\fb(\ell^2(\Gamma))$.

As we shall see this is the key property that is required to make Haagerup's method work in 
a greater generality.

Following Haagerup \cite[Def. 1.6]{Haagerup:Inventiones} we say that a function $\phi:\Gamma\longrightarrow \bbc$ is a multiplier of $C^*_r(\Gamma ) $ if and only if there exists a unique bounded operator 
$M_\phi : C^*_r(\Gamma) \ra C^*_r(\Gamma)$ such that 
\begin{equation}\label{EQN1}
M_\phi \lambda(\gamma) = \phi(\gamma) \lambda (\gamma).
\end{equation}
for all $\gamma \in \Gamma$. This condition can be written equivalently as: 
\begin{equation}\label{EQN2}
M_\phi\lambda(f) = \lambda (\phi\cdot f). 
\end{equation}

An important situation in which such operators arise is given by the following lemma, which 
is a generalisation of \cite[Lemma 1.7]{Haagerup:Inventiones}.

\begin{lemma} \label{multiplier}
Let $\Gamma$ be a discrete group equipped with a length function $\ell$. Assume that $(\Gamma,\ell)$
satisfies the rapid decay inequality for  given $C,s > 0$. 

Let $\phi$ be any function on $\Gamma$ such that
$$
K = \sup_{\gamma\in \Gamma}| \phi(\gamma)| (1 + \ell(\gamma))^s < \infty. 
$$

Then $\phi$ is  a multiplier of  $C^*_r(\Gamma)$ and 
$
\norm{M_\phi} \leq CK. 
$

In particular this holds for any element $f\in \bbc\Gamma$ and for any such element $M_f$ has finite rank.
\end{lemma}
 
 \begin{proof}
 This proof is essentially identical to the proof of Lemma 1.7 in \cite{Haagerup:Inventiones}. 
 For any discrete group $\Gamma$ the characteristic function $\delta_e$ of the identity element 
 $e$ of $\Gamma$ is the identity of the group ring $\bbc \Gamma$. Since $\delta_e$ is a unit 
 vector in $\ell^2(\Gamma) $ we have that for any $f\in \bbc\Gamma$, $\norm{\lambda(f)} 
 \geq \norm{\lambda(f)(\delta_e)}_2 = \norm{f*\delta_e}_2 = \norm{f}_2$. 
 
Then for any 
 $f\in \bbc \Gamma$, the pointwise product  $\phi\cdot f$ is also an element of $\bbc \Gamma$, so 
 we can apply the rapid decay inequality to get: 
 \begin{displaymath}
 \begin{split}
 \norm{\lambda(\phi \cdot f)} &  \leq 
 C \sqrt{\sum_{\gamma\in \Gamma} |\phi(\gamma) f(\gamma)|^2 (1+\ell(\gamma))^{2s}} \\
 & \leq C  \sup_{\gamma\in \Gamma}\{| \phi(\gamma)| (1 + \ell(\gamma))^s \}
 \sqrt{\sum_{\gamma\in \Gamma} |f(\gamma)|^2 } \\
 & = CK\norm{f}_2
 \end{split}
 \end{displaymath}
 
 Putting together the two inequalities we have that 
 $$
 \norm{\lambda(\phi \cdot f)} \leq CK\norm{f}_2 \leq CK \norm{\lambda(f)}. 
 $$
This shows that the map from $\bbc \Gamma$ to $C^*_r(\Gamma)$ which sends $\lambda (f)$ to 
$\lambda(\phi \cdot f)$ is continuous and so extends to a unique map 
$M_\phi : C^*_r(\Gamma) \ra C^*_r(\Gamma)$ with the property that $M_\phi \lambda (f)
= \lambda (\phi \cdot f)$. 

It is also clear that $\norm{M_\phi} \leq CK$. Finally it is clear that if $\phi$ has finite support then $M_\phi$ has finite rank.

 \end{proof}
 
 We are now ready to prove the main technical statement of this note. 
 \begin{theorem}\label{main}
 Let $\Gamma$ be a discrete group with a conditionally negative length function $\ell$, which 
 satisfies the  property (RD) for $C, s > 0$. Then  there exists a net $\{ \phi _\alpha\}$ of functions on $\Gamma$ with finite support such that 
 \begin{enumerate}
 \item For each $\alpha$,  
 $\norm{M_{\phi_\alpha}} \leq 1$; 
 \item $\norm{M_{\phi_\alpha} (x) - x } \to 0$ for all $x\in C^*_r(\Gamma)$. 
 \end{enumerate}
 \end{theorem}
\begin{proof}
 Since the length function  $\ell$ is conditionally negative, it follows from  Schoenberg's lemma that for any $r > 0$ the function $\phi_r(\gamma) = e^{-r \ell(\gamma)}$ is  of positive type. Thus by 
\cite[Lemma 1.1]{Haagerup:Inventiones} (see also \cite[Lemma 3.2 and 3.5]{HV}), for every $r$ there exists a unique completely positive operator
$M_{\phi_r}: C^*_r(\Gamma) \ra C^*_r(\Gamma)$ such that $M_{\phi_r}(\lambda(\gamma))
= \phi_r(\gamma) \lambda (\gamma)$ for all $\gamma \in \Gamma$ and 
$\norm{M_{\phi_r}} = \phi_r(e) = 1$. 

Let us now define a family $\phi_{r,n}$ of finitely supported functions on $\Gamma$ by truncating 
the functions $\phi_r$ to balls of radius $n$ with respect to the length function $\ell$. For every
$\gamma \in \Gamma $ we  put: 
$$
\phi_{r,n} (\gamma) = \begin{cases}
e^{-r\ell(\gamma)}, & \text{if}\;   \ell(\gamma) \leq n \\
0,  & \text{otherwise}.
\end{cases}
$$

Since $e^{-x}(1+x)^s \to 0$ for any positive $s$ and $x\to \infty$, $\sup_{\gamma \in \Gamma}
|\phi_r(\gamma)|(1+\ell(\gamma))^s < \infty$. If we denote this finite number by $K$, then 
clearly $\sup_{\gamma \in \Gamma}
|\phi_{r,n}(\gamma)|(1+\ell(\gamma))^s \leq K$. Thus, for every $r$ and $n$, these functions 
are multipliers of  $C^*_r(\Gamma)$, and the corresponding operators $M_{\phi_r}$ and 
$M_{\phi_{r,n}}$ have norms bounded by $CK$. Since the functions 
$\phi_{r,n}$ have finite support,  the corresponding  operators $M_{\phi_{r,n}}$ are of finite rank.

On the other hand, since
$$
(\phi_r - \phi_{r,n})(\gamma) = \begin{cases} 
0, & \ell(\gamma) \leq n\\
e^{-r \ell(\gamma)}, & \ell(\gamma) > n
\end{cases}
$$
 we have that 
 \begin{displaymath}
 \begin{split}
 \sup_{\gamma \in \Gamma} | (\phi_r  - \phi_{r,n}) & (\gamma)|  (1+\ell(\gamma))^s \\
&  =  \sup_{\ell(\gamma) > n} | (\phi_r - \phi_{r,n})(\gamma)| (1+\ell(\gamma))^s \\
 & \leq K_n < \infty
 \end{split}
 \end{displaymath}
 where $K_n \to 0$ as $n\to \infty$. Thus these functions are multipliers of  $C^*_r(\Gamma)$ 
  and the corresponding operators $M_{\phi_r - \phi_{r,n}}$ are such that 
 $\norm{M_{\phi_r - \phi_{r,n}}} \leq CK_n \to 0$, as $n\to \infty$. 
 
  
Since 
  $$
  \norm{M_{\phi_r} - M_{\phi_{r,n}}} = \norm{M_{\phi_r - \phi_{r,n}}}
  $$
 we have $ \norm{M_{\phi_r} - M_{\phi_{r,n}}} \to 0$ as $n\to \infty$. This implies that 
  $\norm{M_{\phi_{r,n}}} \to \norm{M_{\phi_r}} = \phi_r(e) = 1$. 
  
  To get the correct bound on the norm of these operators we introduce scaled functions: 
  $$
  \rho_{r,n} = \frac{1}{\norm{M_{\phi_{r,n}}}} \phi_{r,n}. 
  $$
 The algebraic identity satisfied by the multipliers, as stated in (\ref{EQN2}),  guarantees that on 
 $\lambda(\bbc \Gamma)$ we have the following identity
 \begin{equation}\label{EQN3}
 M_{\rho_{r,n}} =  \frac{1}{\norm{M_{\phi_{r,n}}}}M_{ \phi_{r,n}}. 
 \end{equation}
 
We now want to show that each operator $M_{\rho_{r,n}}$ is a finite rank contraction
on $C^*_r(\Gamma)$ and  that the strong operator closure of the family 
$\{ M_{\rho_{r,n}} \}$ contains the 
identity map $\text{id}: C^*_r(\Gamma) \ra C^*_r(\Gamma)$. This means that for every 
positive $\epsilon$ there exists an  operator $M_{\rho_{r,n}}$ such that 
$$
\norm{ M_{\rho_{r,n}}x - x} < \epsilon
$$
for all $x\in C^*_r(\Gamma)$. 

First,  a simple use of the triangle inequality leads to the following argument. 
  \begin{equation}\label{EQN3.1}
  \begin{split}
  \norm{M_{\rho_{r,n}} - M_{\phi_r}} & \leq \norm{M_{\rho_{r,n}} - M_{\phi_{r,n}}} + 
  \norm{M_{\phi_{r,n}} - M_{\phi_r}} \\
  & = \norm{(1 - 1/\norm{M_{\phi_{r,n}}}) M_{\phi_{r,n}}}  + \norm{M_{\phi_{r,n}} - M_{\phi_r}} \\
  & \to 0 \qquad \text{as}\quad n\to \infty
  \end{split}
  \end{equation}
  
 Let $x\in C^*_r(\Gamma)$. Then $x$ is a limit of a sequence of elements $x_m \in \lambda(\bbc \Gamma)$ so that
 $
 \norm{M_{\rho_{r,n}}(x)} = \lim_{m\to\infty} \norm{M_{\rho_{r,n}}(x_m)}, 
 $
and equation (\ref{EQN3}) implies that   $\norm{M_{\rho_{r,n}}(x_m)} 
= \norm{
                \frac{1}{\norm{M_{\phi_{r,n}}}} M_{\phi_{r,n}}(x_m)
                }$. 
                
                This leads to the following estimate: 
\begin{equation}\label{EQN3.5}
\begin{split}
\norm{M_{\rho_{r,n}}(x)} & = \lim_{m\to \infty} \norm{\frac{1}{\norm{M_{\phi_{r,n}} }}M_{\phi_{r,n}} (x_m)}\\
&\leq \lim_{m\to \infty} \norm{\frac{1}{\norm{M_{\phi_{r,n}} }}M_{\phi_{r,n}} }\norm{x_m} =\lim_{m\to\infty} \norm{x_m} = \norm{x}.
\end{split}
\end{equation}
%
%


It follows that $\norm{M_{\rho_{r,n}}} \leq 1$.

Finally, it is clear that for any $\gamma \in \Gamma$, $e^{-r\ell(\gamma)}\to 1$ as $r \to 0$. 
Thus for any  $x =   \sum_{\gamma\in \Gamma} \mu_\gamma \lambda (\gamma) \in \bbc \Gamma$ we 
have 
$$
M_{\phi_r}(x) = \sum \mu_\gamma \phi_r(\gamma) \lambda (\gamma)
$$
so that 
\begin{equation}\label{EQN4}
\begin{split}
\lim_{r\to 0} M_{\phi_r}(x)  & = \lim_{r\to 0} \sum \mu_\gamma \phi_r(\gamma)\lambda(\gamma) \\
& = \sum \mu_\gamma (\lim_{r \to 0} \phi_r(\gamma)) \lambda(\gamma) \\
& = \sum \mu_\gamma \lambda(\gamma) = x
\end{split} 
\end{equation}

Since any $x\in C^*_r(\Gamma)$ can be approximated by a sequence $x_m\in \lambda (\bbc \Gamma)$ 
we have
\begin{displaymath}
\begin{split}
\norm{M_{\phi_r}(x) - x } & \leq \norm{M_{\phi_r} (x) - M_{\phi_r}(x_m)} \\
& + \norm{M_{\phi_{r}} (x_m) - x_m} + \norm{x_m - x} \\
\end{split}
\end{displaymath}

Given  that $\norm{M_{\phi_r}}\leq 1$ for all $r > 0$, $\norm{M_{\phi_r} (x) - M_{\phi_r}(x_m)}\leq \norm{x-x_m}<  \epsilon / 3$ for all large enough $n$ and independently of $r$. Thus the sum of the 
first and third term of this sum can be made smaller than $(2/3)\epsilon$, for all $r >0$, and independently of $m$. 
Now equation (\ref{EQN4}) shows  that, as $r\to 0$, $M_{\phi_r}(x_m)$ tends to $x_m$ so the middle term 
will be smaller than $\epsilon/3$ for all sufficiently small $r$. 
Thus, for all sufficiently small $r>0$, $\norm{M_{\phi_r}(x) - x} < \epsilon$ and so
$$
\norm{M_{\phi_r}(x) - x} \to 0
$$
as $r \to 0$ for all $x\in C^*_r(\Gamma)$. 

Let $\epsilon > 0$. Then it follows from (\ref{EQN3.1}) that  for every $r > 0 $ and all sufficiently large 
$n$, $\norm{M_{\rho_{r,n}} - M_{\phi_r}} < \epsilon/2$. Secondly, as we have just shown, for all sufficiently small $r$, 
$\norm{M_{\phi_r}(x) - x} < \epsilon/2$. Given that 
$$
\norm{M_{\rho_{r,n}}x - x} \leq \norm{M_{\rho_{r,n}}x  - M_{\phi_r}x}  + \norm{M_{\phi_r}(x) - x}$$

 for every $x\in C^*_r(\Gamma)$, the norm on the left hand side 
can be made smaller than $\epsilon$ by taking a sufficiently large $n$ and a sufficiently small $r> 0$. 

This means that the strong closure of the family $\fm= \{ M_{\rho_{r,n}}\}$ of finite rank contractions 
contains the identity map on the algebra $C^*_r(\Gamma)$. This implies that there exists 
a net of finitely supported functions $\phi_\alpha$ with corresponding finite rank contractions 
$M_{\phi_\alpha} \in \fm$  such that $\norm{M_{\phi_\alpha}x - x}\to 0$. This concludes the proof. 
\end{proof}


As a corollary we obtain the main result of this note. 
 \begin{theorem}
 Let $\Gamma$ be a discrete group satisfying the rapid decay property with respect to a length function $\ell$
 which is conditionally negative. Then the reduced $C^*$-algebra $C^*_r(\Gamma)$ has the metric 
 approximation property. 
 \end{theorem}
 
Now according to Niblo and Reeves \cite{NR} given a group acting on a CAT(0) cube complex we obtain a conditionally negative kernel on the group which gives rise to a conditionally negative length function. By results of Chatterji and Ruane \cite{CR} the group will have the rapid decay property with respect to this this length function provided that  the action is properly discontinuous, stabilisers are uniformly bounded and the cube complex has finite dimension. Hence we obtain:

 \begin{corollary}
 Groups acting properly discontinuously on a finite dimensional CAT(0) cube complex with uniformly bounded stabilisers have the metric 
 approximation property. 
 \end{corollary}

This class of examples includes free groups, finitely generated Coxeter groups \cite{NRCox}, and finitely generated  right angled Artin groups for which the Salvetti complex is a CAT(0) cube complex. A rich class of interesting examples is furnished by Wise, \cite{W}, in which it is shown that many small cancellation groups  act properly and co-compactly on CAT(0) cube complexes. The examples include  every finitely presented group satisfying the  B(4)-T(4) small cancellation condition and all those  word-hyperbolic groups satisfying the B(6) condition.

Another class of examples where the main theorem applies is furnished by groups acting co-compactly and properly discontinuously on real or complex hyperbolic space. According to a result of Faraut and Harzallah \cite{Crofton} the natural metrics on these hyperbolic spaces are conditionally negative and they give rise to conditionally negative length functions on the groups. See \cite{Robertson} for a discussion and generalisation of this fact. The fact that these metrics satisfy rapid decay for the group was established by Jolissaint in \cite{Jolissaint}.

Finally we remark that  the net $\phi_\alpha$ of Theorem 
\ref{main} provides an approximate identity for the Fourier algebra $A(\Gamma)$ of the group $\Gamma$ which is bounded in the multiplier norm. The proof of \cite[Theorem 2.1]{Haagerup:Inventiones} carries over verbatim to the present situation. This implies, as in 
\cite[Corollary 2.2]{Haagerup:Inventiones}, that  if a group $\Gamma$ satisfies the (RD) property 
with respect to a conditionally negative length function then its Fourier algebra $A(\Gamma)$ has 
the metric approximation property.



\begin{thebibliography}{999}
\bibitem{CR} I.~Chatterji, K.~Ruane, Some geometric groups with rapid decay, \emph{ArXiv preprint math.GR/0310356}, http://xxx.lanl.gov/abs/math.GR/0310356. 
  \bibitem{Haagerup:Inventiones} U.~Haagerup, An example of a non-nuclear $C^*$-algebra which has
  the metric approximation property, Inventiones Math. 50 (1979), 279-293. 
  \bibitem{Crofton}J.~Faraut and K.~Harzallah 
Distances hilbertiennes invariantes sur un espace homogne.Ann. Inst. Fourier (Grenoble) 24 (1974), no. 3, xiv, 171--217.
  \bibitem{HV} P.~de la Harpe, A.~Valette, La propri\'{e}t\'{e} (T) de Kazhdan 
  pour les groupes localement compacts, Asterisque 175 (1989), Soc. Math\'{e}matique de France. 
  \bibitem{Jolissaint} P. Jolissaint. Rapidly decreasing functions in reduced C*-algebras of groups. Trans. Amer. Math. Soc. 317 (1990), 167Ð196. 
\bibitem{NR}G.A.Niblo, L.D.Reeves. The geometry of cube
complexes and the complexity of their fundamental groups,Topology,
Vol. 37,No 3,pp 621-633,1998.
\bibitem{NRCox}G.A.Niblo and L.D.Reeves, Coxeter groups act on CAT(0) cube complexes, Journal of Group Theory, 6, (2003), pp 309-413.
\bibitem{Robertson}G.~Robertson,
Crofton formulae and geodesic distance in hyperbolic spaces. J. Lie Theory 8 (1998), no. 1, 163--172.
\bibitem{W}Daniel T. Wise, Cubulating Small Cancellation Groups, Preprint, http://www.gidon.com/dani/tl.cgi?athe=pspapers/SmallCanCube.ps.  \end{thebibliography}
 \end{document}